\documentclass[11pt,reqno,twoside]{amsart}

\usepackage{amsmath, amssymb, amscd, amsthm, amscd}

\theoremstyle{plain}
\newtheorem{thm}{Theorem}[section]

\newtheorem{lem}[thm]{Lemma}
\newtheorem{prop}[thm]{Proposition}
\newtheorem{Def}[thm]{Definition}

\theoremstyle{remark}
\newtheorem{rem}[thm]{Remark}
\newcounter{sspar}[subsection]
\renewcommand\thesspar{(\thesubsection.\arabic{sspar})}
    {\par\ \newline
     \vskip-\baselineskip\vskip.1truecm
     \noindent\refstepcounter{sspar}
     \noindent\textbf{\thesspar} \ignorespaces}
    {\vskip-\baselineskip
    \ignorespaces}
    {\refstepcounter{sspar}
     \textup{\textbf{\thesspar}} \ignorespaces}
    {\vskip-\baselineskip
    \ignorespaces}
\newcommand{\R}{\mathbb R}
\newcommand{\N}{\mathbb N}
\newcommand{\C}{\mathbb C}
\newcommand{\Z}{\mathbb Z}

\newcommand{\al}{\alpha}
\newcommand{\be}{\beta}

\newcommand{\Ga}{\Gamma}
\newcommand{\de}{\delta}
\newcommand{\De}{\Delta}
\newcommand{\eps}{\varepsilon}

\newcommand{\Te}{\Theta}
\newcommand{\thet}{\vartheta}
\newcommand{\la}{\lambda}

\newcommand{\ph}{\varphi}
\newcommand{\Ups}{\Upsilon}

\newcommand{\Om}{\Omega}
\newcommand{\lt}{\ell^2}
\newcommand{\inprod}[2]{\langle #1,#2 \rangle}
\newcommand{\tensor}{\otimes}
\newcommand{\dirint}{\sideset{}{^\oplus}\int\limits_{0}^{\infty}}

\newcommand{\F}[5]{\,_{#1}F_{#2} \left( \genfrac{.}{.}{0pt}{}{#3}{#4}
\ ;#5 \right)}
\newcommand{\hf}{\frac{1}{2}}
\newcommand{\su}{\mathfrak{su}}

\newcommand{\Conf}{{}_1F_1}
\newcommand{\Res}[1]{\underset{#1}{\mathrm{Res}}}
\newcommand{\D}{\mathcal D}
\newcommand{\vect}[2]{ \begin{pmatrix} #1 \\ #2 \end{pmatrix} }

\numberwithin{equation}{section}

\begin{document}

\date{\today}
\title{Laguerre functions and representations of $\su(1,1)$}
\author{Wolter Groenevelt}
\address{
Technische Universiteit Delft, ITS-TWA \\
Postbus 5031, 2600 GA Delft, The Netherlands}
\email{W.G.M.Groenevelt@its.tudelft.nl}
\dedicatory{Dedicated to Tom Koornwinder on the occasion of his 60th birthday.}

\begin{abstract}
Spectral analysis of a certain doubly infinite Jacobi operator leads to orthogonality relations for confluent hypergeometric functions, which are called Laguerre functions. This doubly infinite Jacobi operator corresponds to the action of a parabolic element of the Lie algebra $\su(1,1)$. 
The Clebsch-Gordan coefficients for the tensor product representation of a positive and a negative discrete series representation of $\su(1,1)$ are determined for the parabolic bases. They turn out to be multiples of Jacobi functions.  From the interpretation of Laguerre polynomials and functions as overlap coefficients, we obtain a product formula for the Laguerre polynomials, given by a discontinuous integral over Laguerre functions, Jacobi functions and continuous dual Hahn polynomials.
\end{abstract}

\maketitle
\section{Introduction}
Many special functions of hypergeometric type have an interpretation in representation theory of Lie groups and Lie algebras, see for example the book by Vilenkin and Klimyk \cite{VK}. In this paper we consider the three dimensional Lie algebra $\su(1,1)$, generated by $H$, $B$ and $C$. 
Elements of $\su(1,1)$ are either elliptic, parabolic or hyperbolic elements, which correspond to the three conjugacy classes of the Lie group $SU(1,1)$. The self-adjoint element $X_a=-aH+B-C$, $a \in \R$, is an elliptic element for $|a|>1$, a parabolic element for $|a|=1$, and a hyperbolic element for $|a|<1$. In \cite{KJ1} Koelink and Van der Jeugt consider the action of $X_a$ in tensor products of positive discrete series representations. This leads to convolution identities for several hypergeometric orthogonal polynomials. The idea to look at elements of the form of $X_a$ is due to Granovskii and Zhedanov \cite{GZ}. The action of $X_a$ in the tensor product of a positive and a negative discrete series representation is considered in \cite{GK} for the elliptic case, and in \cite{GKR} for the hyperbolic case. In this paper we investigate the remaining, parabolic, case. In \cite{WG1} the quantum version of $X_a$ is studied. The Lie algebra $\su(1,1)$ is replaced by the quantized universal enveloping algebra $\mathcal U_q\big(\su(1,1)\big)$ and $X_a$ is replaced by a twisted primitive element. It turns out that in $\mathcal U_q\big(\su(1,1)\big)$ the three cases are all the same.

There are four classes of irreducible unitary representations of $\su(1,1)$, the positive and negative discrete series, the principal unitary series and the complementary series. The tensor product of a positive and a negative discrete series representation decomposes into a direct integral over the principal unitary series. Discrete terms can occur, and these terms correspond to one complementary series, or a finite number of discrete series. The Clebsch-Gordan coefficients for the standard bases are multiples of continuous dual Hahn polynomials.

We consider the element $X = -H+B-C$, which is a parabolic element. In the discrete series $X$ acts on the standard (elliptic) basis as a Jacobi operator, which corresponds to the three term recurrence relation for Laguerre polynomials.  In the principal unitary series and the complementary series $X$ acts on the standard basis as a doubly infinity Jacobi operator, which corresponds to the recurrence relation for Laguerre functions. So the Laguerre polynomials and functions appear as overlap coefficients between the (generalized) eigenvectors of $X$ and the standard basis vectors. Using the differential equation for the Laguerre polynomials, we realize the generators $H$, $B$ and $C$ in the discrete series as differential operators acting on polynomials. In these realizations the action of the Casimir operator can be identified with the hypergeometric differential equation, which leads to Jacobi functions as Clebsch-Gordan coefficients for parabolic basis vectors. This gives a product formula for Laguerre polynomials, which has a similar structure as the discontinuous integral for Bessel functions of Weber and Schafheitlin.

This paper is organized as follows. In  \S\ref{app:Laguerre} we consider a certain doubly infinity Jacobi operator, which corresponds to the action of $X$ in the principle unitary series. Spectral analysis leads to orthogonality relations for Laguerre functions. This section is based on \cite{MR} by Masson and Repka.

In \S\ref{sec3} we turn to representations of the Lie algebra $\su(1,1)$. We introduce the orthogonal polynomials and functions that we need in \S\ref{ssec:functions}, and we give some of their properties. In \S\ref{ssec:su11} we introduce the Lie algebra $\su(1,1)$ and give the irreducible unitary representations. In \S\ref{ssec:P} we diagonalize the element $X$ in the various representations, and we give generalized eigenvectors. In \S\ref{ssec:CGCs} the generators $H$, $B$ and $C$ are realized as differential operators. Then the Casimir operator in the tensor product can be identified with the hypergeometric differential operator, and this gives Jacobi functions as Clebsch-Gordan coefficients. As a result we obtain a product formula for Laguerre polynomials, which involves Jacobi functions, Laguerre functions and continuous dual Hahn polynomials. 

\emph{Notations.}
If $d\mu(x)$ is a positive measure, we use the notation $d\mu^\hf(x)$ for the positive measure with the property 
\[
d(\mu^\hf \times \mu^\hf)(x,x) = d\mu(x).
\]

The hypergeometric series is defined by
\[
\F{p}{q}{a_1, \ldots, a_{p}} {b_1, \ldots, b_q}{z} = \sum_{n=0}^\infty \frac{ (a_1)_n \ldots (a_{p})_n }{ (b_1)_n \ldots (b_q)_n } \frac{z^n}{n!},
\]
where $(a)_n$ denotes the Pochhammer symbol, defined by
\[
(a)_n = \frac{\Ga(a+n)}{\Ga(a)}=a(a+1)(a+2) \ldots (a+n-1), \qquad n \in \Z_{\geq 0}.
\]
For the confluent hypergeometric function we use the notation
\[
\Conf(a;b;z) = \F{1}{1}{a}{b}{z},
\]
and the second solution of the confluent hypergeometric differential equation is defined by 
\begin{equation} \label{def:U}
U(a;b;z) = \frac{ \Ga(1-b)}{\Ga(a-b+1)} \Conf(a;b;z) + \frac{ \Ga(b-1)}{\Ga(a)}z^{1-b} \Conf(a-b+1;2-b;z),
\end{equation}
see \cite[(1.3.1)]{Sl}. This is a many valued functions of $z$, and we take as its principal branch that which lies in the complex plane cut along the negative real axis $(-\infty,0]$.

\section{Laguerre functions} \label{app:Laguerre}
In this section we determine the spectral measure of a certain doubly infinite Jacobi operator. This operator is obtained from the action of the self-adjoint element $X$ of the Lie algebra $\su(1,1)$ in the principal unitary series representation, see \S\ref{ssec:P}. The eigenfunctions which are needed to describe the spectral measure, are called Laguerre functions. See \cite{MR} or \cite{Koe} for doubly infinite Jacobi operators. The calculation of the eigenfunctions and the Wronskian is obtained from \cite{MR}, but we repeat the calculations here briefly. \\

The doubly infinite Jacobi operator $L:\ell^2(\Z) \rightarrow \ell^2(\Z)$ is defined by
\begin{equation} \label{def:L}
L e_k = a_k e_{k+1} + b_{k} e_k + a_{k-1} e_{k-1},
\end{equation}
where $\{e_k\}_{k \in \Z}$ is the standard orthonormal basis of $\ell^2(\Z)$ and
\[
\begin{split}
a_k = a_k(\rho, \eps)& =| (k+\eps+\hf+i\rho)|, \\
b_k = b_k(\rho,\eps)&= 2(k+\eps),
\end{split}
\]
where $\rho \geq 0$, $\eps \in [0,1)$ and $(\rho,\eps) \neq (0,\hf)$. 

\begin{rem} \label{rem:symm}
There exists a symmetry for the parameters of $L$.
Let us denote $L=L(\rho,\eps)$. The unitary operator $U : e_k \mapsto (-1)^k e_{-k}$ intertwines $L(\rho,\eps)$ with $-L(-\rho,-\eps)$. So if $f(z;\eps,\rho,k)$ is a solution to the eigenvalue equation $Lf=zf$, then $(-1)^k f(-z;-\eps,-\rho,-k)$ is another solution to the same eigenvalue equation.
\end{rem}

Observe that $L$ is an unbounded, symmetric operator. The domain $\mathcal{D}$ of $L$ is the dense subspace of finite linear combinations of the basis vectors $e_k$. 
We define for a function $f = \sum_{k=-\infty}^\infty f_k e_k \in \ell^2(\Z)$ 
\[
L^* f = \sum_{k=-\infty}^\infty (a_kf_{k+1} + b_k f_k + a_{k-1}f_{k-1}) e_k.
\]
on its domain
\[
\D^* = \left\{ f \in \ell^2(\Z)\ | \ L^* f \in \ell^2(\Z) \right\},
\]
then $(L^*, \D^*)$ is the adjoint of $(L, \D)$. Note that $L^*|_\D=L$.

Solutions to $Lv=-zv$ can be given in terms of confluent hypergeometric functions (see \cite{Sl}).
\begin{prop} \label{prop:eigenf L}
The following functions are solutions to $Lv=-zv$:
\[
\begin{split}
s_k(z;\rho,\eps) &= (-1)^k \frac{|\Ga(k+\eps+\hf+i\rho)| }{ \Ga(k+\eps+\hf-i\rho)  }  \, \Conf(k+\eps+\hf+i\rho;1+2i\rho;z),\\
t_k(z;\rho,\eps) &=  \frac{ |\Ga(\hf-k-\eps-i\rho) |}{ \Ga(\hf-k-\eps+i\rho) } \,
\Conf(\hf-k-\eps-i\rho;1-2i\rho;-z),\\
u_k(z;\rho,\eps) &= (-1)^k |\Ga(k+\eps+\hf+i\rho)|\,
 U(k+\eps+\hf+i\rho;1+2i\rho;z), \quad z \not\in (-\infty,0]\\
v_k(z;\rho,\eps) &= |\Ga(\hf-k-\eps-i\rho)|\,
U(\hf-k-\eps-i\rho;1-2i\rho;-z), \quad z\not\in [0,\infty).
\end{split}
\]
\end{prop}
\begin{proof}
The first solution $s_k$ follows from \cite[(2.2.1)]{Sl}
\[
(b-a)\,\Conf(a-1;b;z)+(2a-b+z)\,\Conf(a;b;z) - a\,\Conf(a+1;b;z)=0.
\]
The second solution $t_k$ follows from the first using the symmetry relation for the parameters, cf.~ Remark \ref{rem:symm}. In the same way we find from \cite[(2.2.8)]{Sl} 
\[
U(a-1;b;z)-(2a-b+z)U(a;b;z)+a(a-b+1)U(a+1;b;z)=0,
\]
that $u_k$ and $v_k$ are solutions to $Lv=-zv$.
\end{proof}
The solution space to $Lv=-zv$ is two-dimensional, since for a fixed $n\in \Z$, $v$ is completely determined by the initial values $v_{n-1}$ and $v_n$.
So the eigenfunctions given in Proposition \ref{prop:eigenf L} can be expanded in terms of each other.
\begin{prop} \label{prop:connec}
We have the connection formulas
\[
\begin{split}
u_k(z) =& A(z) s_k(z) + B(z) t_k(z), \quad z \not\in (-\infty,0] \\
v_k(z) =& C(z) s_k(z) + D(z) t_k(z), \quad z \not\in [0,\infty),
\end{split}
\]
where
\[
\begin{split}
A(z) =& \Ga(-2i\rho),\\
B(z) =& z^{-2i\rho} e^{z} \Ga(2i\rho)  \frac{ \sin \pi(\eps+\hf+i\rho) }{| \sin \pi(\eps+\hf-i\rho)|},\\
C(z) =& (-z)^{2i\rho} e^{-z} \Ga(-2i\rho)  \frac{ \sin \pi(\eps+\hf+i\rho) }{ |\sin \pi(\eps+\hf-i\rho)|},\\
D(z) =& \Ga(2i\rho).
\end{split}
\]
Or equivalently, for $z \in \C\setminus\R$,
\[
\begin{split}
s_k(z) &= E(z)u_k(z) + F(z)v_k(z), \\
t_k(z) &= G(z)u_k(z) + H(z)v_k(z),
\end{split}
\]
where
\[
\begin{split}
E(z) &= \frac{1}{ \pi} \sin \pi(\eps+\hf-i\rho) \Ga(1+2i\rho) e^{i\pi \xi(\eps+\hf+i\rho)},\\
F(z) &= -\frac{1}{\pi}\, | \sin \pi (\eps+\hf+i\rho)|\,
\Ga(1+2i\rho) e^z z^{-2i\rho}e^{i\pi \xi(\eps+\hf+i\rho)},\\
G(z) &= \frac{1}{\pi}\,| \sin \pi (\eps+\hf+i\rho) |\,
\Ga(1-2i\rho) e^{-z} (-z)^{2i\rho}e^{i\pi \xi(\eps+\hf+i\rho)},\\
H(z) &=-\frac{1}{ \pi} \sin \pi(\eps+\hf-i\rho) \Ga(1-2i\rho) e^{i\pi \xi(\eps+\hf+i\rho)} ,
\end{split}
\]
where $\xi = \mathrm{sgn}\Im (z)$.
\end{prop}
\begin{proof}
The first connection formula follows from \eqref{def:U}, the reflection formula for the $\Ga$-function, and Kummer's transformation: $\Conf(a;b;z) = e^z \Conf(b-a;b;-z)$. The second connection formula can be derived from the first using the symmetry for the parameters, see Remark \ref{rem:symm}.

The other two connection formulas follow from \cite[(1.9.1),(1.4.10)]{Sl} or they can be derived from the first two.
\end{proof}

\begin{Def} \label{def:Wr}
For two functions $f(z) = \sum_{k=-\infty}^\infty f_k(z) e_k$ and $g(z) = \sum_{k=-\infty}^\infty g_k(z) e_k$, the Wronskian is defined by
\[
[f(z),g(z)]_k = a_k(f_k(z) g_{k+1}(z) - f_{k+1}(z)g_k(z)).
\]
\end{Def}
If $f(z)$ and $g(z)$ are solutions to the eigenvalue equation $Lv=zv$, the Wronskian $[f(z),g(z)]_k$ is is independent of $k$, so the Wronskian can be found by taking the limit $k \rightarrow \pm \infty$. Moreover, $f(z)$ and $g(z)$ are linearly independent solutions if and only if $[f(z),g(z)]\neq0$.

\begin{lem} \label{lem:asymp}
For $0<|\arg(z)|<\pi$ and $k \rightarrow \infty$
\[
\begin{split}
u_k(z) &= (-1)^k \sqrt{\pi}\, e^{\hf z-\sqrt{4(k+\eps)z}} z^{-\frac{1}{4}-i\rho} k^{-\frac{1}{4}} \Big(1+ \mathcal O(k^{-\hf}) \Big),\\
v_{-k}(z) &= \sqrt{\pi}\,e^{-\hf z-\sqrt{-4(k-\eps)z}} (-z)^{-\frac{1}{4}+i\rho} k^{-\frac{1}{4}}  \Big(1+ \mathcal O(k^{-\hf}) \Big),\\
s_k(z) &= \frac{(-1)^k}{2\sqrt{\pi}}\, e^{\hf z+\sqrt{4(k+\eps)z}} z^{-\frac{1}{4}-i\rho} \Ga(1+2i\rho) k^{-\frac{1}{4}}\Big(1+ \mathcal O(k^{-\hf}) \Big),\\
t_{-k}(z) &= \frac{1}{2\sqrt{\pi}}\,e^{-\hf z+\sqrt{-4(k-\eps)z}} (-z)^{-\frac{1}{4}+i\rho} \Ga(1-2i\rho) k^{-\frac{1}{4}} \Big(1+ \mathcal O(k^{-\hf}) \Big).
\end{split}
\]
\end{lem}
\begin{proof}
This follows from the asymptotic behaviour for $|y| \rightarrow \infty$ of the modified Bessel functions
\[
I_\nu(y) = \frac{ e^y }{ \sqrt{2 \pi y}}\left(1+\mathcal O\Big(\frac{1}{|y|}\Big)\right), \quad K_\nu(y) = \sqrt{ \frac{ \pi}{2y}} e^{-y}\left(1+\mathcal O\Big(\frac{1}{|y|}\Big)\right),\qquad |\arg y|<\frac{\pi}{2},
\]
and the asymptotic expansions for the confluent hypergeometric functions in terms of modified Bessel functions \cite[(4.6.42),(4.6.43)]{Sl}.
\end{proof}

For $z \in \C \setminus \R$ we introduce the spaces
\[
\begin{split}
S^+_z &= \left\{ f(z)=\sum_{k=-\infty}^\infty f_k(z) e_k\ |\ L^*f(z) = -zf(z) \text{ and } \sum_{k=0}^\infty |f_k(z)|^2 < \infty \right\},\\
S^-_z& = \left\{ f(z)=\sum_{k=-\infty}^\infty f_k(z) e_k\ |\ L^*f(z) = -zf(z) \text{ and } \sum_{k=-\infty}^{-1} |f_k(z)|^2 < \infty \right\}.
\end{split}
\]
Then the deficiency space for $L$ is $S^+_z \cap S^-_z$. Note that $\mathrm{dim} \ S_z^\pm  \leq 2$, and in case $\mathrm{dim} \ S_z^\pm  = 2$, we have $S_z^+ = S_z^-$, since the solution space to $Lf = -zf$ is two dimensional.

Next we put $\phi_z = z^{i\rho} u(z)$, $z \not\in (\infty,0]$, and $\Phi_z = (-z)^{-i\rho}v(z)$, $z \not\in[0,\infty)$. From the transformation $U(a;b;z) = z^{1-b}U(a-b+1;2-b;z)$ it follows that
\begin{equation} \label{eq:compl}
\begin{split}
u_l(\overline{z}) &= \overline{z}^{-2i\rho} \overline{u_l(z)}, \quad z \not\in (-\infty,0],\\
v_k(\overline{z})& = (-\overline{z})^{2i\rho} \overline{v_k(z)}, \quad z \not\in [0,\infty).
\end{split}
\end{equation}
So we have $\overline{(\phi_z)_k} = (\phi_{\overline{z}})_k$ and $\overline{(\Phi_z)_k} =
(\Phi_{\overline{z}})_k$, and in particular $(\phi_{x})_k \in \R$ for $x>0$ and $(\Phi_{x})_k \in \R$ for $x<0$. Note that $(\phi_{x})_k$ and $(\Phi_{x})_k$ are even in $\rho$. 

We calculate the Wronskian $[\phi_z, \Phi_z]$. 
From Lemma \ref{lem:asymp} and $a_k=k+\mathcal O(1)$, for $k \rightarrow \infty$, we find for $ z \not\in (-\infty,0]$
\[
[s(z),u(z)] = \lim_{k\rightarrow \infty} \hf e^z z^{-\hf-2i\rho} \Ga(1+2i\rho) k^{\hf}\Big( e^{\sqrt{4(k+\eps+1)z}-\sqrt{4(k+\eps)z}}- e^{\sqrt{4(k+\eps)z}-\sqrt{4(k+\eps+1)z}}\Big).
\]
And since
\[
\Big( e^{\sqrt{4(k+\eps+1)z}-\sqrt{4(k+\eps)z}}- e^{\sqrt{4(k+\eps)z}-\sqrt{4(k+\eps+1)z}}\Big) = 2 \sqrt{\frac{z}{k}}\left(1+ \mathcal O(k^{-\hf})\right),\qquad k\rightarrow \infty,
\]
we obtain
\[
[s(z),u(z)]= e^{z} z^{-2i\rho} \Ga(2i\rho+1), \qquad z \not\in (-\infty,0].
\]
Then we find from the connection formulas of Proposition \ref{prop:connec}
\[
[s(z),u(z)] = F(z)[v(z),u(z)] ,\qquad z \in \C \setminus \R,
\]
and this gives
\begin{equation} \label{eq:Wr}
[\phi_z,\Phi_z] = z^{i\rho}(-z)^{-i\rho}[u(z),v(z)]=\frac{-\pi e^{-i \pi\xi(\eps+\hf)}}
{ | \sin \pi(\eps+\hf+i\rho)|}, \qquad 0<|\arg(z)|<\pi.
\end{equation}
So we find that $\phi_z$ and $\Phi_z$ are linearly independent.
\begin{prop} \label{S+S-}
For $0<|\arg(z)|<\pi$, we have $S_z^+= \mathrm{span} \left\{ \phi_z \right\}$, and $S_z^-= \mathrm{span}\left\{ \Phi_z \right\}$, and $L$ is essentially self-adjoint.
\end{prop}
\begin{proof}
From Lemma \ref{lem:asymp} we see that $\phi_z = z^{i\rho}u(z) \in S_z^+$ and $\Phi_z=(-z)^{-i\rho}v(z) \in S_z^-$, for $0<|\arg(z)|< \pi$. Masson and Repka prove in \cite[Thm.2.1]{MR} that the deficiency indices of $L$ are obtained by adding the deficiency indices of the two Jacobi operators $J^\pm$ obtained by restricting $L$ to $\ell^2(\Z_{\geq 0})$ (setting $a_{-1}=0$) and to $\ell^2(-\N)$ (setting $a_0=0$). Since $\sum_{k=0}^\infty 1/a_k$ and $\sum_{k=-\infty}^{-1} 1/a_k$ are divergent, \cite[Ch.VII, Thm.1.3]{Be} proves that $J^\pm$ have deficiency indices $(0,0)$, and hence so has $L$. So $\mathrm{dim}\ S^\pm_z = 1$, and the proposition follows. 
\end{proof}
We use the Stieltjes-Perron inversion formula, see \cite[\S XII.4]{DS}, to calculate the spectral measure;
\[
E_{f,g}\Big((a,b)\Big) = \lim_{\eta \downarrow 0} \lim_{\de \downarrow 0} \frac{1}{2\pi i}\int_{a+\eta}^{b-\eta}
\langle G(x+i\de)f,g \rangle - \langle G(x-i\de)f,g \rangle dx.
\]
In this case the resolvent $G(z)$ can be calculated explicitly by
\begin{equation} \label{def:resol}
\inprod{ G(z)f}{g} = \frac{1}{[\phi_z,\Phi_z]}
\sum_{k\leq l}(\Phi_z)_k (\phi_z)_l (f_l \overline{g}_k + f_k
\overline{g}_l)(1-\hf \de_{lk}).
\end{equation}

\begin{prop} \label{prop:sp meas}
The spectral measure for the operator $-L$ defined by \eqref{def:L}, is described by the following integral, for $\mathcal B \subseteq \R$ a Borel set,
\[
\begin{split}
\inprod{E(\mathcal B)f}{g} =& \frac{1}{\pi^2}\int_{\mathcal B \cap (-\infty,0)} e^{x}|\sin \pi (\eps+\hf+i\rho) |^2 \inprod{f}{v(x)} \inprod{v(x)}{g} dx \\
+&\frac{1}{\pi^2}\int_{\mathcal B \cap (0,\infty)} e^{-x}|\sin \pi (\eps+\hf+i\rho) |^2 \inprod{f}{u(x)} \inprod{u(x)}{g} dx,
\end{split}
\]
where $f, g \in \lt(\Z)$ and with the notation of Proposition \ref{prop:eigenf L}
\[
u(x) = \sum_{k=-\infty}^\infty u_k(x) e_k, \quad v(x) = \sum_{k=-\infty}^\infty v_k(x) e_k.
\]
Moreover, $0$ is not contained in the point spectrum of $L$. 
\end{prop}
\begin{proof}
We define 
\begin{equation} \label{def:De(x)}
\De(x)=\lim_{\de\downarrow 0}\ \frac{ (\Phi_{x+i\de})_k (\phi_{x+i\de})_l }{ [\phi_{x+i\de},\Phi_{x+i\de}]} -
\frac{ (\Phi_{x-i\de})_k (\phi_{x-i\de})_l }{ [\phi_{x-i\de},\Phi_{x-i\de}]},  
\end{equation}
then we have, using \eqref{eq:Wr} and $\xi= \mathrm{sgn}\Im(z)$,
\[
\De(x) = -\lim_{\de\downarrow 0}\ \frac{1}{\pi}\, | \sin \pi(\eps+\hf+i\rho) |\,\Big( e^{i\pi(\eps+\hf)} (\Phi_{x+i\de})_k (\phi_{x+i\de})_l - e^{-i\pi(\eps+\hf)} (\Phi_{x-i\de})_k
(\phi_{x-i\de})_l\Big).
\]
For $x<0$ we have $\lim_{\de \downarrow 0} (\Phi_{x+i\de})_k = \lim_{\de \downarrow 0} (\Phi_{x-i\de})_k= (\Phi_x)_k$, and from the connection formulas we find
\[
\begin{split}
\lim_{\de\downarrow 0}\ &e^{i\pi(\eps+\hf)} (\phi_{x+i\de})_l- e^{-i\pi(\eps+\hf)} (\phi_{x-i\de})_l \\
=&\, \lim_{\de\downarrow 0}  \Big(e^{i\pi(\eps+\hf)} (x+i\de)^{i\rho}A(x+i\de)s_l(x+i\de) - e^{-i\pi(\eps+\hf)} (x-i\de)^{i\rho}A(x-i\de)s_l(x-i\de) \Big)\\
&+ \lim_{\de\downarrow 0}\Big(e^{i\pi(\eps+\hf)} (x+i\de)^{i\rho}B(x+i\de)t_l(x+i\de) - e^{-i\pi(\eps+\hf)} (x-i\de)^{i\rho}B(x-i\de)t_l(x-i\de) \Big)\\
=&\, 2i\sin \pi(\eps+\hf+i\rho)(-x)^{i\rho}\Ga(-2i\rho)s_l(x)
 + 2i e^{x} (-x)^{-i\rho} \Ga(2i\rho)\ |\sin \pi(\eps+\hf+i\rho)|\, t_l(x)\\
=&\, 2i e^{x} (-x)^{-i\rho} |\sin \pi(\eps+\hf+i\rho)|\, v_l(x).
\end{split}
\]
Here we used
\[
\lim_{\de \downarrow 0} (-y\pm i \de)^a = e^{\pm i\pi a}y^a, \quad y >0.
\]
For $x>0$ we use the symmetry for the parameters, cf.~ Remark \ref{rem:symm}. So we find
\begin{equation} \label{eq:De(x)}
\De(x)=
\begin{cases}
\displaystyle -\frac{2i}{\pi}  e^x (-x)^{-2i\rho}|\sin\pi(\eps+\hf+i\rho)|^2 v_k(x)v_l(x),&x<0,\\
&\\
\displaystyle -\frac{2i}{\pi} e^{-x}x^{2i\rho}|\sin \pi (\eps+\hf+i\rho) |^2 u_k(x)u_l(x), &x>0.
\end{cases}
\end{equation}
Both expressions are clearly symmetric in $k$ and $l$, so the sum in \eqref{def:resol} can be antisymmetrized using \eqref{eq:compl}. Now, if $0$ is not contained in the point spectrum of $L$, the result follows from the Stieltjes-Perron inversion formula and \eqref{eq:compl}.

To show that $0$ is not an element of the point spectrum of $L$, we show that $\mathrm{ker} L = \{0\}$. First we calculate the Wronskian $[s(0),t(0)]$, using Definition \ref{def:Wr} with $k=0$. A straightforward calculation gives
\[
[s(0),t(0)] = 2i\rho \frac{ \sin \pi (\eps+\hf-i\rho) }{ | \sin \pi( \eps+\hf+i\rho) |},
\]
hence $s(0)$ and $t(0)$ are linearly independent. So if $f \in \mathrm{ker}L$, $f \neq 0$, then $f$ is a linear combination of $s(0)$ and $t(0)$. But $s(0), t(0) \not\in \ell^2(\Z)$, since $|s_k(0)|=1$ and $|t_k(0)|=1$, and therefore $\mathrm{ker}L = \{0\}$.
\end{proof}
\begin{rem}
The result of Proposition \ref{prop:sp meas} remains valid if $i\rho$ is replaced by $\la+\hf$ where $\la \in (-\hf, -\eps)$ and $\eps \in [0,\hf)$, or $\la \in (-\hf, \eps-1)$ and $\eps \in (\hf,1)$. In this case operator $L$ is obtained from the action of the self-adjoint element $X$ in the complementary series representation of $\su(1,1)$, see \S\ref{ssec:P}.
\end{rem}

Let us define the Laguerre functions $\psi_n(x;\rho,\eps)$, $n \in \Z$, by
\[
\psi_n(x;\rho,\eps) =
\begin{cases}
v_n(x;\rho,\eps), & x<0,\\ 
u_n(x;\rho,\eps), & x>0,
\end{cases}
\]
and we define the weight function $w(x;\rho,\eps)$ by
\[
w(x;\rho,\eps) = 
\frac{1}{\pi^2} |\sin \pi (\eps+\hf+i\rho) |^2 e^{-|x|}.
\]
From Proposition \ref{prop:sp meas} we find the following.
\begin{thm} \label{thm:Lag func}
For $\rho \geq 0$, $\eps \in [0,1)$ and $(\rho,\eps) \neq (0,\hf)$, the Laguerre functions $\psi_n(x;\rho,\eps)$ form an orthonormal basis of $L^2(\R, w(x;\rho,\eps)dx)$.
\end{thm}
\begin{proof}
The orthonormality of the Laguerre functions follows from Proposition \ref{prop:sp meas} by replacing $f$ and $g$ by standard orthonormal basis vectors $e_n$ and $e_m$, $n,m \in \Z$, and using $\mathcal B = \R$. Completeness of the Laguerre functions follows from the uniqueness of the spectral measure.
\end{proof}

We also want to define the Laguerre functions for $x=0$. In order to find the natural definition in $x=0$ we calculate $\De(0)$, where $\De(x)$ is defined by \eqref{def:De(x)}. Using the connection coefficients from Proposition \ref{prop:connec} we find
\[
\begin{split}
\De&(0)=-\lim_{\de\downarrow 0}\ \frac{1}{\pi}\, |\sin \pi (\eps+\hf+i\rho)|
\Big( e^{i\pi(\eps+\hf)} (\Phi_{i\de})_k (\phi_{i\de})_l - e^{-i\pi(\eps+\hf)} (\Phi_{-i\de})_k
(\phi_{-i\de})_l\Big)\\
&=-\lim_{\de\downarrow 0}\ \frac{1}{\pi}\, |\sin \pi (\eps+\hf+i\rho)|
\Big( e^{i\pi(\eps+\hf+i\rho)}  \big[ C(i\de) s_k(i\de)+D(i\de)t_k(i\de) \big]\, \big[ A(i\de) s_l(i\de)+B(i\de) t_l(i\de) \big]\\
&\qquad \qquad  -e^{-i\pi(\eps+\hf+i\rho)}  \big[ C(-i\de) s_k(-i\de)+D(-i\de)t_k(-i\de) \big]\, \big[ A(-i\de) s_l(-i\de)+B(-i\de) t_l(-i\de) \big]\Big).
\end{split}
\]
To compute this limit, we use 
\[
\begin{split}
\lim_{\de\downarrow 0} s_k(i\de) &= \lim_{\de\downarrow 0} s_k(-i\de) =s_k(0),\\
\lim_{\de\downarrow 0} t_k(i\de) &= \lim_{\de\downarrow 0} t_k(-i\de) =t_k(0),\\
\lim_{\de\downarrow 0} A(i\de) D(i\de) &= \lim_{\de\downarrow 0} A(-i\de) D(-i\de) = \left|\Ga(2i\rho) \right|^2,\\
\lim_{\de\downarrow 0} B(i\de) C(i\de)& =  e^{-i\pi\, 2i\rho} \left|\Ga(2i\rho) \right|^2 \frac{ \sin \pi (\eps+\hf+i\rho) }{ \sin \pi (\eps+\hf-i\rho) },\\
\lim_{\de\downarrow 0} B(-i\de) C(-i\de) &=  e^{i\pi\, 2i\rho} \left|\Ga(2i\rho) \right|^2 \frac{ \sin \pi (\eps+\hf+i\rho) }{ \sin \pi (\eps+\hf-i\rho) },
\end{split}
\]
then we find
\[
\De(0) =- \frac{2i}{\pi} \sin\pi (\eps+\hf+i\rho)\,| \sin \pi (\eps+\hf+i\rho) | \left|\Ga(2i\rho) \right|^2 \big( t_k(0) s_l(0) + s_k(0)t_l(0)\big).
\]
From Euler's reflection formula for the $\Ga$-function we obtain 
\[
t_k(0) =  \frac{ \sin \pi( \eps+\hf -i\rho) }{| \sin \pi( \eps+\hf 
+i\rho) |}\ \overline{ s_k(0)}, 
\]
and this gives
\[
\begin{split}
\De(0) &=-\frac{2i}{\pi}\,| \sin \pi (\eps+\hf+i\rho)\,\Ga(2i\rho) |^2 \big(t_k(0) \overline{t_l(0)}+\overline{t_k(0)}t_l(0)   \big)\\
&=-\frac{2i}{\pi}\,| \sin \pi (\eps+\hf+i\rho)\, \Ga(2i\rho) |^2 \vect{\overline{t_k(0)}}{t_k(0)}^* \vect{\overline{t_l(0)}}{t_l(0)}.
\end{split}
\] 
Comparing this result with \eqref{eq:De(x)}, we see that for $x=0$ the Laguerre function can be defined by
\[
\psi_n(0;\rho,\eps)=|\Ga(2i\rho)|\vect{\overline{t_k(0) }}{  t_k(0)} .
\]

\section{Clebsch-Gordan coefficients for parabolic basis vectors of $\su(1,1)$} \label{sec3}

\subsection{Orthogonal polynomials and functions} \label{ssec:functions}
The Wilson polynomials, see Wilson \cite{Wil} or \cite[\S3.8]{AAR},
are polynomials on top of the Askey-scheme of hypergeometric polynomials,
see Koekoek and Swarttouw \cite{KS}.
The continuous dual Hahn polynomials are a three-parameter subclass of
the Wilson polynomials, and are defined by
\begin{equation} \label{def:cont dHahn}
s_n(y;a,b,c) = (a+b)_n (a+c)_n\F{3}{2}{-n,a+ix,a-ix}{a+b,a+c}{1},
\quad x^2=y.
\end{equation}
For real parameters $a$, $b$, $c$, with $a+b$, $a+c$, $b+c$ positive,
the continuous dual Hahn polynomials are orthogonal with respect to a
positive measure, supported on a subset of $\R$. The orthonormal
continuous dual Hahn polynomials are defined by
\[
S_n(y;a,b,c) = \frac{(-1)^n s_n(y;a,b,c)} {\sqrt{n!(a+b)_n
(a+c)_n (b+c)_n}}
\]
By Kummer's transformation, see e.g. \cite[Cor. 3.3.5]{AAR}, the
polynomials $s_n$ and $S_n$ are symmetric in $a$, $b$ and $c$.
Without loss of generality we assume that $a$ is the smallest of the
real parameters $a$, $b$ and $c$. Let $d\mu(\cdot;a,b,c)$ be the
measure defined by
\begin{align*}
\int_{\R} f(y) d&\mu(y;a,b,c)  =
 \frac{1}{2\pi}\int_0^\infty \left|
\frac{\Gamma(a+ix)\Gamma(b+ix)\Gamma(c+ix)}{\Gamma(2ix)} \right|^2
\frac{f(x^2)}{\Gamma(a+b) \Gamma(a+c) \Gamma(b+c)}\,d x\\ + &
\frac{\Gamma(b-a) \Gamma(c-a)}{\Gamma(-2a) \Gamma(b+c)} \sum_{k=0}^K
(-1)^k \frac{(2a)_k (a+1)_k (a+b)_k (a+c)_k} {(a)_k (a-b+1)_k (a-c+1)_k
k!} f(-(a+k)^2),
\end{align*}
where $K$ is the largest non-negative integer such that $a+K<0$. In
particular, the measure $d\mu(\cdot;a,b,c)$ is absolutely continuous
if $a\geq0$. The measure is positive under the conditions $a+b>0$,
$a+c>0$ and $b+c>0$. Then the polynomials $S_n(y;a,b,c)$ are
orthonormal with respect to the measure $d\mu(y;a,b,c)$.\\

The Laguerre polynomials are defined by 
\begin{equation} \label{def:Lag pol}
L^{(\al)}_n(x) = \frac{(\al+1)_n}{n!} \Conf(-n;\al+1;x).
\end{equation}
The orthonormal Laguerre polynomials 
\[
l_n^{(\al)}(x) = \sqrt{ \frac{n!} {(\al+1)_n }}\,L^{(\al)}_n(x) .
\]
are orthonormal on $[0,\infty)$ with respect to the weight function
\[
w^{(\al)}(x)=\frac{ x^\al e^{-x} }{ \Ga(\al+1)}.
\]
They satisfy the three-term recurrence relation
\[
xl_n^{(\al)}(x) = - \sqrt{ (n+1)(\al+n+1)} l_{n+1}^{(\al)}(x) + (2n+\al+1) l_n^{(\al)}(x) - \sqrt{ n(n+\al) } l_{n-1}^{(\al)}(x),
\]
and the differential equation
\[
xy''(x) + (\al+1-x)y'(x)+ny(x) =0,\quad y(x) = l_n^{(\al)}(x).
\]

The Jacobi functions, see \cite{Koo}, are defined by
\begin{equation} \label{def:Jac funct}
\ph_\la^{(\al,\be)}(x) = \F{2}{1}{\hf(\al+\be+1-i\la), \hf(\al+\be+1+i\la)}{\al+1}{-x}.
\end{equation}
Here we use the unique analytic continuation to $\C \setminus [1,\infty)$ of the hypergeometric function.
The Jacobi functions are eigenfunctions of the hypergeometric differential operator
\[
-x(1+x)\frac{ d^2}{dx^2}-[\al+1+(\al+\be+2)x]\frac{d}{dx}  
\]
for eigenvalue $\frac{1}{4}[(\al+\be+1)^2+\la^2]$. Spectral analysis of the hypergeometric differential operator leads to a unitary integral transform called the Jacobi-function transform. The Jacobi-function transform is given by
\begin{equation} \label{def: Jacobi trans}
\begin{cases}
\displaystyle(\mathcal F f)(\la)= \int_0^\infty f(x)\ph_\la^{(\al,\be)}(x) \De_{\al,\be}(x) dx \\ \\
\displaystyle f(x) = \frac{1}{2\pi}\int (\mathcal F f)(\la)\ph_\la^{(\al,\be)}(x) d\nu(\la)
\end{cases}
\end{equation}
where $\al>-1$, $\be \in \R$, $\De_{\al,\be}(x)= 2^{2\al+2\be+1} x^{\al}(1+x)^{\be}$, and $d\nu(\la)$ is the measure given by
\[
\begin{split}
\frac{1}{2\pi}\int g(\la) d\nu(\la) &= \frac{1}{2\pi}\int_0^\infty g(\la) |c_{\al,\be}(\la)|^{-2} d\la -i \sum_{\la \in \mathcal D} g(\la) \Res{\mu=\la} \big(c_{\al,\be}(\mu) c_{\al,\be}(-\mu) \big)^{-1}, \\
c_{\al,\be}(\la) &= \frac{2^{\al+\be+1-i\la} \Ga(\al+1) \Ga(i\la) }{ \Ga\big(\hf(\al+\be+1+i\la)\big) \Ga\big(\hf(\al-\be+1+i\la)\big) },\\
\mathcal D &= \Big\{ i(|\be|-\al-1-2j)\ | \ j \in \Z_{\geq 0}, |\be|-\al-1-2j>0 \Big\}.
\end{split}
\]
Observe that the measure $d\nu(\la)$ is absolutely continuous if $|\be|\leq \al+1$.\\

\subsection{The Lie algebra $\boldsymbol{\su(1,1)}$} \label{ssec:su11}
The Lie algebra $\su(1,1)$ is a three dimensional Lie algebra, generated by $H$, $B$ and $C$ satisfying the commutation relations
\begin{equation} \label{commutation rel}
[H,B]=2B, \quad [H,C]=-2C, \quad [B,C]=H.
\end{equation}
There is a $\ast$-structure by $H^*=H$ and $B^*=-C$. The Casimir operator $\Omega$ is a central element of $\mathcal U\big(\su(1,1)\big)$, and $\Om$ is given by
\begin{equation} \label{Casimir}
 \Omega = -\frac{1}{4}(H^2+2H+4CB) .
\end{equation}

There are four classes of irreducible unitary representations of
$\su(1,1)$, see \cite[\S6.4]{VK}:

The positive discrete series
representations $\pi_k^+$ are representations labelled by $k>0$. The
representation space is $\lt(\Z_{\geq 0})$ with orthonormal basis
$\{e_n\}_{n\in \Z_{\geq 0}}$. The action is given by
\begin{equation} \label{pos}
\begin{split}
\pi_k^+(H)\, e_n =&\ 2(k+n)\, e_n,  \\ 
\pi_k^+(B)\, e_n =&\ \sqrt{(n+1)(2k+n)}\, e_{n+1},  \\ 
\pi_k^+(C)\, e_n =&\ -\sqrt{n(2k+n-1)}\, e_{n-1},  \\ 
\pi_k^+(\Omega)\, e_n =&\ k(1-k)\, e_n.
\end{split}
\end{equation}

The negative discrete series representations $\pi_{k}^-$ are labelled by
$k>0$. The representation space is $\lt(\Z_{\geq 0})$ with orthonormal
basis $\{e_n\}_{n \in \Z_{\geq 0}}$. The action is given by
\begin{equation} \label{neg}
\begin{split}
\pi_{k}^-(H)\, e_n =&\ -2(k+n)\, e_n, \\ 
\pi_{k}^-(B)\, e_n =&\ -\sqrt{n(2k+n-1)}\, e_{n-1},  \\ 
\pi_{k}^-(C)\, e_n =&\ \sqrt{(n+1)(2k+n)}\, e_{n+1},  \\ 
\pi_{k}^-(\Omega)\, e_n =&\ k(1-k)\, e_n.
\end{split}
\end{equation}

The principal series representations $\pi^{\rho,\eps}$ are labelled by
$\eps \in [0,1)$ and $\rho \geq 0$, where $(\rho,\eps) \neq
(0,\frac{1}{2})$. The representation space is $\lt(\Z)$ with
orthonormal basis $\{e_n\}_{n \in \Z}$. The action is given by
\begin{equation} \label{princ}
\begin{split}
\pi^{\rho, \eps}(H)\, e_n =&\ 2(\eps+n)\, e_n, \\ 
\pi^{\rho, \eps}(B)\, e_n =&\ |(n+\eps+\frac{1}{2}+i\rho)|\, e_{n+1},\\ 
\pi^{\rho, \eps}(C)\, e_n =&\ -|(n+\eps-\frac{1}{2}+i\rho)|\, e_{n-1},
\\ 
\pi^{\rho, \eps}(\Omega)\, e_n =&\ (\rho^2+\frac{1}{4})\, e_n.
\end{split}
\end{equation}
For $(\rho,\eps)= (0,\hf)$ the representation $\pi^{0,\hf}$ splits into a direct sum of a positive and a negative discrete series representation: $\pi^{0,\hf}= \pi^+_\hf \oplus \pi^-_\hf$. The representation space splits into two invariant subspaces: $\{e_n \, |\, n<0 \} \oplus \{ e_n\, | \, n \geq 0 \}$.

The complementary series representations $\pi^{\lambda,\eps}$ are
labelled by $\eps$ and $\lambda$, where $\eps \in [0,\frac{1}{2})$ and
$\lambda \in (-\frac{1}{2},-\eps)$  or $\eps \in (\frac{1}{2},1)$ and
$\lambda \in (-\frac{1}{2},\eps-1)$. The representation space is
$\lt(\Z)$ with orthonormal basis $\{e_n\}_{n \in \Z}$. The action is
given by
\begin{equation} \label{comp}
\begin{split}
\pi^{\lambda,\eps}(H)\,e_n =&\ 2(\eps+n)\, e_n, \\ 
\pi^{\lambda,\eps}(B)\, e_n =&\ \sqrt{(n+\eps+1+\lambda) (n+\eps-\lambda)} \, e_{n+1}, \\
\pi^{\lambda,\eps}(C)\, e_n =&\ -\sqrt{(n+\eps+\lambda) (n+\eps-\lambda-1)}\, e_{n-1},  \\ 
\pi^{\lambda,\eps}(\Omega)\, e_n =&\ -\lambda(1+\lambda)\, e_n.
\end{split}
\end{equation}
Note that formally for $\la = -\frac{1}{2}+i\rho$ the actions in the principal
series and in the complementary series are the same.\\

We remark that the operators \eqref{pos}-\eqref{comp} are unbounded, with domain the set of finite linear combinations of the basis vectors. The representations are $*$-representations in the sense of Schm\"udgen \cite[Ch.8]{Sch}.\\

The decomposition of the tensor product of a positive and a
negative discrete series representation of $\su(1,1)$ is determined in
\cite[Thm.2.2]{GK}, see also \cite[\S8.7.7]{VK} for the group $SU(1,1)$.
\begin{thm} \label{thm:decomp}
For $k_1 \leq k_2$ the decomposition of the tensor product of positive
and negative discrete series representations of $\su(1,1)$ is
\begin{align*}
\pi_{k_1}^+ \tensor\pi_{k_2}^- &\cong \dirint \pi^{\rho,\eps} d \rho, &
k_1-k_2 \geq -\frac{1}{2}, k_1+k_2 \geq \frac{1}{2}, \\ \pi_{k_1}^+
\tensor\pi_{k_2}^- &\cong \dirint \pi^{\rho,\eps} d \rho \oplus
\pi^{\la, \eps}, & k_1+k_2<\frac{1}{2}, \\ \pi_{k_1}^+
\tensor\pi_{k_2}^- &\cong \dirint \pi^{\rho,\eps} d \rho \oplus
\bigoplus_{\substack{j \in \Z_{\geq 0}\\ k_2-k_1-\frac{1}{2}-j>0}}
\pi_{k_2-k_1-j}^-, & k_1-k_2
< -\frac{1}{2},
\end{align*}
where $\eps = k_1-k_2+L$, $L$ is the unique integer such that $\eps \in
[0,1)$, and $\la = -k_1-k_2$. Further, under the identification above,
\begin{equation} \label{decomp vec}
e_{n_1} \tensor e_{n_2} = (-1)^{n_2} \int_\R S_n(y;n_1-n_2)
e_{n_1-n_2-L}\, d\mu^\hf(y;n_1-n_2),
\end{equation}
where $n= \min\{n_1,n_2\}$, $S_n(y;p)$ is an orthonormal continuous dual Hahn
polynomial,
\[
S_n(y;p) =
\begin{cases}
S_n(y;k_1-k_2+\frac{1}{2},k_1+k_2-\frac{1}{2},k_2-k_1-p+\frac{1}{2}),&
p \leq 0,\\
S_n(y;k_2-k_1+\frac{1}{2},k_1+k_2-\frac{1}{2},k_1-k_2+p+\frac{1}{2}),&
p\geq0,
\end{cases}
\]
and
\[
d \mu(y;p) =
\begin{cases}
d
\mu(y;k_1-k_2+\frac{1}{2},k_1+k_2-\frac{1}{2},k_2-k_1-p+\frac{1}{2}),&
p \leq 0,\\ d \mu
(y;k_2-k_1+\frac{1}{2},k_1+k_2-\frac{1}{2},k_1-k_2+p+\frac{1}{2}),&
p\geq0.
\end{cases}
\]
\end{thm}
The inversion of
\eqref{decomp vec} can be given explicitly, e.g. for an element
\[
f \tensor e_{r-L} = \int_0^\infty f(x) e_{r-L} dx \in L^2(0,\infty)
\tensor \lt(\Z) \cong \dirint \lt(\Z) dx
\]
in the representation space of the direct integral representation, we have
\begin{equation} \label{inv decomp}
f \tensor e_{r-L} =
\begin{cases}
\displaystyle \sum_{p=0}^\infty (-1)^{p-r} \left[ \int_{\R} S_p(y;r)
f(y) d\mu^\hf(y;r) \right] e_p \tensor e_{p-r}, & r \leq 0,\\ \displaystyle
\sum_{p=0}^\infty (-1)^{p} \left[\int_{\R} S_p(y;r) f(y) d \mu^\hf(y;r)
\right] e_{p+r} \tensor e_p, & r \geq 0.
\end{cases}
\end{equation}
For the discrete components in Theorem \ref{thm:decomp} we can replace $f$
by a Dirac delta function at the appropriate points of the discrete
mass of $d \mu(\cdot;r)$. In the following subsections we assume that discrete terms do not occur in the tensor product decomposition. From the calculations it is clear how to extend the results to the general case. \\

\subsection{Parabolic basisvectors} \label{ssec:P}

We consider the self-adjoint element
\[
X = -H+B-C \in \su(1,1),
\]
which is a parabolic element. We determine the spectral decomposition of $X$ in the various representations. We also give (generalized) eigenvectors of $X$.  This is done in the same way as in \cite{KJ1}, using (doubly infinite) Jacobi operators. First we consider $X$ in the discrete series. The action of $X$ can be identified with the three-term recurrence relation for the Laguerre
polynomials.
\begin{prop} \label{prop:+-}
The operators $\Te^{\pm}$ defined by
\[
\begin{split}
\Te^\pm: \ell^2(\Z_{\geq 0}) &\rightarrow L^2\big([0,\infty), w^{(2k-1)}(x)dx) \\
e_n &\mapsto l_n^{(2k-1)}(\cdot),
\end{split}
\] 
are unitary and intertwine $\pi^\pm_k(X)$ with $M_{\mp x}$.
\end{prop}
Here $M$ denotes the multiplication operator: $M_fg(x)= f(x) g(x)$.

In terms of generalized eigenvectors, Proposition \ref{prop:+-} states that
\[
v^\pm(x) = \sum_{n=0}^\infty l_n^{(2k-1)}(x) e_n, \qquad x \in [0,\infty),
\]
is a generalized eigenvector of $\pi^\pm_k(X)$ for eigenvalue $\mp x$. These eigenvectors can be considered as parabolic basis vectors for $\su(1,1)$.

Next we consider $X$ in the principal unitary series. We find that $\pi^{\rho,\eps}(X)$ extends to a doubly infinite Jacobi operator which corresponds to the recurrence relation for the Laguerre functions. The spectral analysis of $\pi^{\rho,\eps}(X)$ is carried out in \S\ref{app:Laguerre}. 
\begin{prop} \label{prop:princ}
The operator $\Te^{\rho,\eps}$ defined by
\[
\begin{split}
\Te^{\rho,\eps} : \ell^2(\Z) &\rightarrow L^2\big(\R,w(x;\rho,\eps)dx\big)\\
e_n &\mapsto 
(-1)^n \psi_n(\cdot;\rho,\eps)
\end{split}
\]
is unitary and intertwines $\pi^{\rho,\eps}(X)$ with $M_{x}$.
\end{prop} 
So, for $x \in \R$,
\[
v^{\rho,\eps}(x) =
 \sum_{n=-\infty}^\infty (-1)^n \psi_n(x;\rho,\eps)
e_n, \\
\]
is a generalized eigenvector of $\pi^{\rho,\eps}(X)$ for
eigenvalue $x$. Note that $v^{\rho,\eps}(0)$ is two dimensional, since $\psi_n(0)$ is given by a two dimensional vector.

Next we consider the action of $X$ in the tensor product. Recall that $\De(Y) = 1 \tensor Y +Y \tensor 1$ for $Y \in \su(1,1)$. Then we find from Proposition \ref{prop:+-} the following.
\begin{prop} \label{prop:+-1}
The operator $\Ups$ defined by
\[
\begin{split}
\Ups: \ell^2(\Z_{\geq 0})\tensor  \ell^2(\Z_{\geq 0})& \rightarrow L^2\big([0,\infty)\times [0,\infty), w^{(2k_1-1)}(x_1)w^{(2k_2-1)}(x_2)dx_1 dx_2 \big)\\
e_{n_1} \tensor e_{n_2} & \mapsto l_{n_1}^{(2k_1-1)}(x_1)l_{n_2}^{(2k_2-1)}(x_2),
\end{split}
\] 
is unitary and intertwines $\pi^+_{k_1} \tensor \pi^-_{k_2}\big(\De(X)\big)$ with $M_{x_2-x_1}$.
\end{prop}
So 
\[
v^+(x_1) \tensor v^-(x_2) = \sum_{n_1,n_2=0}^\infty l_{n_1}^{(2k_1-1)}(x_1) l_{n_2}^{(2k_2-1)}(x_2)\, e_{n_1} \tensor e_{n_2}
\]
is a generalized eigenvector of $\pi^+_{k_1} \tensor \pi^-_{k_2}\big( \De(X)\big)$ for eigenvalue $x_2-x_1$.

\subsection{Clebsch-Gordan coefficients} \label{ssec:CGCs}
We want to determine the Clebsch-Gordan coefficients between the uncoupled eigenvectors $v^+(x_1) \tensor v^-(x_2)$ and the coupled eigenvectors $\int^\oplus v^{\rho,\eps}(x_1-x_2) d\rho$. This comes down to finding the function $g$, for which the operator $\Ups_g$ defined by
\[
\begin{split}
\Ups_g : L^2(0,\infty)
\tensor \lt(\Z) \cong \dirint \lt(\Z) dx &\rightarrow \dirint L^2\left(\R,   w(t;\rho,\eps)dt \right) d\rho,\\
f \tensor e_n &  \mapsto (-1)^n \int_0^\infty f(\rho) g(\rho) \psi_n(t;\rho,\eps)d\rho,
\end{split}
\]
is the same as $\Ups$. From Theorem \ref{thm:decomp} and \eqref{princ} we know that the Clebsch-Gordan coefficient $g$ must be an eigenfunction of the Casimir operator $\Om$ in the tensor product for eigenvalue $\rho^2+\frac{1}{4}$. So first we determine the actions of the generators $H$, $B$ and $C$ on parabolic basis vectors.

We start with a very simple lemma, which is based on the fact that $\mathfrak{sl}(2,\C)$ is semi-simple, so $[\mathfrak{sl}(2,\C), \mathfrak{sl}(2,\C)] = \mathfrak{sl}(2,\C)$.
\begin{lem} \label{lem:HX}
\[
B = \frac{1}{4}[H,X]+\hf X + \hf H, \quad C= \frac{1}{4}[H,X] - \hf X - \hf H.
\]
\end{lem}
\begin{proof}
From the definition of $X$ and the commutation relations \eqref{commutation rel} we find $[H,X]=2B+2C$. This proves the lemma.
\end{proof}
This lemma shows that to find the action of the generators $H$, $B$ and $C$, it is enough to find the action of $H$, since the action of $X$ is known.

\begin{prop} \label{prop:rea+-P}
In the positive discrete series, the generators $H$, $B$, $C$ have a realization as differential operators acting on polynomials:
\[
\begin{split}
\pi^+_k(H) &= -2x \frac{d^2}{dx^2} - 2(2k-x) \frac{d}{dx} +2k, \\
\pi^+_k(B) &= -x \frac{d^2}{dx^2} -2(k-x) \frac{d}{dx} +(2k-x), \\
\pi^+_k(C) &= x \frac{d^2}{dx^2} + 2k \frac{d}{dx}.
\end{split}
\]
In the negative discrete series $H$, $B$, $C$ have a realization as differential operators on polynomials:
\[
\begin{split}
\pi^-_k(H) &= 2x \frac{d^2}{dx^2} +2(2k-x) \frac{d}{dx} -2k, \\
\pi^-_k(B) &= x \frac{d^2}{dx^2} + 2k \frac{d}{dx}, \\
\pi^-_k(C) &= -x \frac{d^2}{dx^2} - 2(k-x) \frac{d}{dx} +(2k-x).
\end{split}
\]
\end{prop}
\begin{proof}
We show that $\Te^+$ intertwines the actions of $H$, $B$ and $C$ given by \eqref{pos}, with the differential operators given in the proposition.

From the differential equation for the Laguerre polynomials, we find for the action of $H$
\[
\begin{split}
\Te^+ \pi^+_k(H) e_n &=  (2n+2k) l_n^{(2k-1)}(x) \\
&=  -2x \frac{d^2}{dx^2}l_n^{(2k-1)}(x) - 2(2k-x)
\frac{d}{dx} l_n^{(2k-1)}(x) + 2k\, l_n^{(2k-1)}(x) \\
&=\Big(-2x \frac{d^2}{dx^2}  -2(2k-x)\frac{d}{dx}  +2k \Big) \Te^+ e_n.
\end{split}
\]
So we have realized $\pi^+_k(H)$ as a differential operator.
By Proposition \ref{prop:+-} $\pi^+_k(X)$ is realized as the multiplication operator $M_{-x}$. A direct calculation shows that
\[
\Te^+ \pi^+_k([H,X])e_n = \Big(4x \frac{d}{dx} + 2(2k-x)\Big) \Te^+ e_n.
\]
Then Lemma \ref{lem:HX} proves the proposition for the positive discrete series. We find the action in the negative discrete series in the same way, or we use the Lie-algebra isomorphism
$\thet$, given by
\[
\thet(H) = -H, \quad \thet(B) = C, \quad \thet(C)=B.
\]
Then $\pi^+_k\big(\thet(Y)\big) = \pi^-_k(Y)$ for $Y \in \su(1,1)$.

A straightforward calculation shows that these operators indeed
satisfy the $\su(1,1)$ commutation relations.
\end{proof}

It is also possible to find the actions of $H$, $B$ and $C$ on the
eigenvectors $v^{\rho,\eps}(x)$. This is done using the differential equation for the Laguerre functions, which follows from the confluent hypergeometric differential equation. We do not need these actions here.\\

Next we want to calculate $\pi^+_{k_1} \tensor \pi^-_{k_2} \big( \De(\Om) \big)$ for these realizations. From \eqref{Casimir} we obtain
\begin{equation} \label{De Om}
\De(\Om) = 1 \tensor \Om + \Om \tensor 1 - \hf H \tensor H - ( C \tensor B + B \tensor C).
\end{equation}
Proposition \ref{prop:rea+-P} shows that $\pi^+_{k_1} \tensor
\pi^-_{k_2} \big( \De(\Om) \big)$ is a differential operator acting on polynomials in two variables. Let $p(x_1)$ and $q(x_2)$ be polynomials in $x_1$, respectively $x_2$, then we find after a long calculation in which many terms cancel,
\begin{equation} \label{real De(OM)P}
\begin{split}
\pi^+_{k_1} \tensor \pi^-_{k_2}& \big( \De(\Om) \big) p(x_1) q(x_2) =\\
-&x_1 x_2 \Big( p''(x_1)q(x_2) + 2p'(x_1) q'(x_2) + p(x_1) q''(x_2) \Big)\\ +&
(2x_1x_2-2k_1x_2-2k_2x_1) \Big( p'(x_1)q(x_2) + p(x_1) q'(x_2) \Big) \\
+& \big(2k_1x_2+2k_2x_1-2k_1k_2-x_1x_2+k_1(1-k_1)+k_2(1-k_2)\big) p(x_1) q(x_2),
\end{split}
\end{equation}
where $p'(x_1)= \frac{d}{dx_1}p(x_1)$ and $q'(x_2) =
\frac{d}{dx_2}q(x_2)$. We show that $\pi^+_{k_1} \tensor
\pi^-_{k_2} \big( \De(\Om) \big)$ can be identified with the hypergeometric
differential operator, and therefore the eigenfunctions are Jacobi
functions. Let $t \in \R$ and define 
\[
P(x,t) = 
\begin{cases}
p(x-t) q(x), & t<0,\\
p(x) q(x), & t=0,\\
p(x) q(x+t), & t>0.
\end{cases}
\]
\begin{prop} \label{prop: CG coef P}
The operator $\Xi$, defined by
\[
\Xi:P(x,t) \mapsto
\begin{cases}
e^{x}\ph_{2\rho}^{(2k_2-1, 2k_1-1)}(-x/t), & t<0,\\
e^{x} x^{\hf-k_1-k_2 \pm i\rho} & t=0,\\
e^{x}\ph_{2\rho}^{(2k_1-1, 2k_2-1)}(x/t), & t>0,
\end{cases}
\]
intertwines $\pi^+_{k_1} \tensor
\pi^-_{k_2} \big( \De(\Om) \big)$ with $M_{\rho^2+1/4}$.
\end{prop}
\begin{proof}
First we assume $t> 0$. Put $x=x_1$, $t=x_2-x_1$ and $e^{x}\ph(x;t) =  p(x)q(x+t)$ in \eqref{real De(OM)P}, then
\[
\begin{split}
\pi^+_{k_1} \tensor \pi^-_{k_2}& \big( \De(\Om) \big)e^x \ph(x;t) = \\&
e^x\Big[-x(x+t) \frac{d^2\ph}{dx^2} - (2k_1(x+t)+2k_2x) \frac{d\ph}{dx}+ \big((k_1(1-k_1)+k_2(1-k_2)-2k_1k_2 \big) \ph \Big].
\end{split}
\]
Comparing this to the differential equation for the Jacobi
functions,  we see that the Jacobi function $\ph_{2\rho}^{(2k_1-1, 2k_2-1)}(x/t)$ is an eigenfunction of this differential operator for eigenvalue $\rho^2+\frac{1}{4}$. And then the intertwining property of $\Xi$ follows. 
The case $t<0$ is proved similarly.

Next we assume $t=0$. We put $x = x_1 =x_2$ and $e^{x}\ph(x) =  p(x)q(x)$ in \eqref{real De(OM)P}, and then from Theorem \ref{thm:decomp} it follows that we must solve the following Euler differential equation:
\[
-x^2 y'' - (2k_1+2k_2)x y' + (k_1+k_2)(1-k_1-k_2) y = (\rho^2+\frac{1}{4})y.
\]
The general solutions to this equation are given by $y = c_1 x^{\hf-k_1-k_2+i\rho} +c_2 x^{\hf-k_1-k_2-i\rho}$.
\end{proof}
For $t \neq 0$ we can identify the spectrum of $\pi^+_{k_1} \tensor
\pi^-_{k_2} \big( \De(\Om) \big)$ with the support of the measure $d\nu(2\rho)$. Note that, naturally, the support is exactly the same as the support of the orthonormality measure for the continuous dual Hahn polynomials given in Theorem \ref{thm:decomp}. 

Since the Jacobi function transform is unitary, there exists a constant $c$ such that $c\, \Xi$ is unitary for $t\neq 0$. So Proposition \ref{prop: CG coef P} determines for $t\neq 0$  the Clebsch-Gordan coefficients for the eigenvectors up to a factor independent of $x$. To determine the exact Clebsch-Gordan decomposition for the eigenvectors, we need to find this factor. For $t=0$ the Laguerre functions are $\C^2$-valued, so the Clebsch-Gordan coefficients for $t=0$ are also $\C^2$-valued. We show that in this case the unitarity of $c\, \Xi$ corresponds to the unitarity of the Mellin transform.

\begin{thm} \label{thm:CG decomp P}
The Clebsch-Gordan coefficients for the parabolic bases are given by 
\[
g(x_1-x_2)=
\begin{cases}
\displaystyle C_-(\rho) e^{x_2} \ph_{2\rho}^{(2k_2-1,2k_1-1)}\left(\frac{x_2}{x_1-x_2} \right),  & x_2-x_1<0,\\
\vect{\overline{C_0(\rho)}\, e^{x_1} x_1^{\hf-k_1-k_2+i\rho}}{ C_0(\rho)\, e^{x_1} x_1^{\hf-k_1-k_2-i\rho} },& x_2-x_1=0,\\
\displaystyle  C_+(\rho) e^{x_1} \ph_{2\rho}^{(2k_1-1,2k_2-1)}\left(\frac{x_1}{x_2-x_1} \right), & x_2-x_1>0,
\end{cases}
\]
where
\[
\begin{split}
C_-(\rho) &= \frac{(-1)^L}{\sqrt{2\pi} } (x_1-x_2)^{k_1+k_2-\hf+i\rho} \sqrt{ \frac{ \Ga(2k_1) } {\Ga(2k_2)}}\left| \frac{  \Ga(k_1+k_2-\hf+i\rho) } {\Ga(k_1-k_2-\hf+i\rho)\Ga(2i\rho) } \right|,\\
C_0(\rho) &= \frac{1}{\sqrt{\pi}}\,  \frac{(-1)^L\,\sqrt{2 \Ga(2k_1) \Ga(2k_2) }\,\Ga(k_1+k_2-\hf+i\rho)}{ |\Ga(k_1+k_2-\hf+i\rho) \Ga(k_1-k_2+\hf+i\rho)|  \Ga(k_2-k_1+\hf-i\rho) },\\
C_+(\rho) &=\frac{1}{\sqrt{2\pi}}(x_2-x_1)^{k_1+k_2-\hf-i\rho} \sqrt{ \frac{ \Ga(2k_2) } {\Ga(2k_1)}} \left| \frac{ \Ga(k_1+k_2-\hf+i\rho) } {\Ga(k_2-k_1+\hf+i\rho)  \Ga(2i\rho) }\right|.
\end{split}
\]
\end{thm}

\begin{proof}
Recall that the Clebsch-Gordan coefficients are the functions $g$ such that $\Ups= \Ups_g$. Put $x=x_1$, $t=x_2-x_1$ and assume $t>0$. Since the Clebsch-Gordan coefficients do not depend on $n_1$ and $n_2$ it is enough to show that for the function $g$ defined in the theorem, we have $\Ups(e_0 \tensor e_0) =\Ups_g\big( \int_0^\infty e_{-L}d\mu^{1/2} \big)$, where $d\mu$ is the orthogonality measure for the continuous dual Hahn polynomials as in Theorem \ref{thm:decomp}. Explicitly, we must prove the following identity
\[
\begin{split}
1 =& \frac{(-1)^L}{ \sqrt{2\pi} } \int_0^\infty e^x\ph_{2\rho}^{(2k_1-1, 2k_2-1)}\left( \frac{x}{t} \right)  \psi_{-L}(t;\rho, k_1-k_2+L)\\ & \times
\frac{ C_+(\rho)}{ \sqrt{ \Ga(2k_1) \Ga(2k_2)}} \left| \frac{ \Ga(k_1+k_2-\hf+i\rho) \Ga(k_2-k_1+\hf+i\rho) \Ga(k_1-k_2+\hf+i\rho) } { \Ga(2i\rho) } \right| d\rho.
\end{split}
\]

We use an integral representation for the second solution of the confluent hypergeometric differential equation, see \cite[(3.2.55)]{Sl},
\[
U(a;b;z) = \frac{z^{a-c}}{\Ga(c)} \int_0^\infty e^{-zy} y^{c-1} \F{2}{1}{a, 1+a-b}{c}{-y}dy, \qquad \Re c>0, \Re z>0,
\]
with parameters given by 
\[
a=k_1-k_2+\hf+i\rho, \quad b=1+2i\rho, \quad c=2k_1, \quad y=\frac{x}{t}, \quad z=t.
\] 
By the definition of the Laguerre function $\psi_n(t)$ for $t>0$, see \S\ref{app:Laguerre}, the definition of a Jacobi function \eqref{def:Jac funct} and Euler's transformation \cite[(2.2.7)]{AAR}, we have
\[
\begin{split}
&\psi_{-L}(t;\rho, k_1-k_2+L) = \\
&(-1)^L 2^{3-4k_1-4k_2} t^{\hf-k_1-k_2+i\rho} \frac{|\Ga(k_1-k_2+\hf +i\rho)|} { \Ga(2k_1) } \int_0^\infty  e^{-x} \ph_{2\rho}^{(2k_1-1, 2k_2-1)}\left( \frac{x}{t} \right) \De_{2k_1-1, 2k_2-1} \left( \frac{x}{t} \right)\frac{dx}{t}.
\end{split}
\] 
Taking the inverse Jacobi transform of this, gives
\[
\begin{split}
1 = 
(-1)^L\frac{t^{k_1+k_2-\hf}}{2\pi}&\int_0^\infty e^{x} \ph_{2\rho}^{(2k_1-1, 2k_2-1)}\left( \frac{x}{t} \right) t^{-i\rho} \psi_{-L}(t;\rho, k_1-k_2+L) \\
&\times \frac{|\Ga(k_1-k_2+\hf+i\rho)|}{\Ga(2k_1) } \left| \frac{  \Ga(k_1+k_2-\hf+i\rho)} {\Ga(2i\rho) } \right|^2
d\rho
\end{split}
\]
This is the desired identity. For $t=x_1-x_2<0$ the theorem is proved similarly.

For  $x=x_1=x_2$ we determine the ($\C^2$-valued) function $g$ for which we have $\Ups(e_{n_1} \tensor e_{n_2}) = \Ups_g \big( (-1)^{n_2}\int_0^\infty S_n\, e_{n_1-n_2-L}d\mu^{1/2} \big)$, see Theorem \ref{thm:decomp}. We start with \cite[(3.7.4)]{Sl};
\[
\int_0^\infty e^{-x} x^{s-1} \Conf(a;b;x) \Conf(c;d;x) dx = \frac{ \Ga(s) \Ga(d) \Ga(d-c-s) }{\Ga(d-c) \Ga(d-s) } \F{3}{2}{ a, s, 1+s-d}{ b, 1+s+c-d}{1}, 
\]
where $\Re(s)>0$ and we assume $a$ to be a non-positive integer. We consider this as a Mellin transform, and then the inverse transform gives, for $k>0$,
\[
 \frac{1}{2\pi i} \int_{k-i\infty}^{k+i\infty} \frac{ \Ga(s) \Ga(d) \Ga(d-c-s) }{ \Ga(d-c) \Ga(d-s) } \F{3}{2}{ a, s, 1+s-d}{ b, 1+s+c-d}{1} x^{-s}ds =  e^{-x} \Conf(a;b;x) \Conf(c;d;x).
\]
We use 
\[
a=-{n_1}, \quad b=2k_1, \quad c=-n_2, \quad d=2k_2, \quad s= k_1+k_2-\hf+i\rho,
\]
and then we transform the terminating $_3F_2$-series using the first formula on page 142 in \cite{AAR}. Next we add to the integral the same integral with $\rho$ replaced by $-\rho$, then we obtain an integral for which the integrand is even in $\rho$. Using Euler's reflection formula and the definitions for the Laguerre polynomials \eqref{def:Lag pol}, the continuous dual Hahn polynomials \eqref{def:cont dHahn}, and the Laguerre functions of argument $0$, we have
\[
\begin{split}
l_{n_1}^{(2k_1-1)}(x)\, l_{n_2}^{(2k_2-1)}(x)& =\\
 \frac{1}{\pi}&\int_0^\infty S_{n_1}(\rho^2;k_1-k_2+\hf, k_1+k_2-\hf, k_2-k_1+n_2-n_1+\hf)\\  
\times&\frac{(-1)^{n_1}\,e^x}{\sqrt{(2k_2)_{n_2-n_1} (n_2-n_1)!}}\,|  \Ga(k_2-k_1+\hf+n_2-n_1+i\rho)| \\
\times& \Bigg(  \frac{ \Ga(k_1+k_2-\hf+i\rho) }{ \Ga(k_2-k_1+\hf-i\rho)} t_{n_1-n_2-L}(0;\rho, k_1-k_2+L)\ x^{\hf-k_1-k_2-i\rho}\\
 &+  \frac{ \Ga(k_1+k_2-\hf-i\rho) }{ \Ga(k_2-k_1+\hf+i\rho)}\overline{ t_{n_1-n_2-L}(0;\rho, k_1-k_2+L)}\ x^{\hf-k_1-k_2+i\rho} \Bigg) d\rho.
\end{split}
\]
Writing the last integral as
\[
(-1)^{n_1+L}\int_0^\infty S_{n_1}(\rho^2;n_1-n_2) \vect{g_1(\rho)}{g_2(\rho)}^* \vect{\overline{t_{n_1-n_2-L}(0)}}{t_{n_1-n_2-L}(0)} d\mu^\hf(\rho;n_1-n_2),
\]
gives the expression for the Clebsch-Gordan coefficients.
\end{proof}

\begin{rem}
The explicit expressions for the Clebsch-Gordan coefficients as $_2F_1$-series can also be found in Basu and Wolf \cite{BW}. The method used in \cite{BW} to compute the Clebsch-Gordan coefficients is different from the method used here.
\end{rem}
Theorem \ref{thm:CG decomp P} gives the following Clebsch-Gordan decomposition for the parabolic basis vectors
\[
v^+(x_1) \tensor v^-(x_2) = 
\begin{cases}
\displaystyle \int_0^\infty C_-(\rho) e^{x_2} \ph_{2\rho}^{(2k_2-1,2k_1-1)}\left(\frac{x_2}{x_1-x_2} \right) v^{\rho,\eps}(x_2-x_1) d\rho, & x_2-x_1<0,\\ \\
\displaystyle \int_0^\infty \vect{\overline{C_0(\rho)}\, e^{x_1} x_1^{\hf-k_1-k_2+i\rho}}{ C_0(\rho)\, e^{x_1} x_1^{\hf-k_1-k_2-i\rho} }^* v^{\rho,\eps}(0)d\rho, &x_2-x_1=0,\\  \\
\displaystyle \int_0^\infty C_-(\rho) e^{x_1} \ph_{2\rho}^{(2k_1-1,2k_2-1)}\left(\frac{x_1}{x_2-x_1} \right) v^{\rho,\eps}(x_2-x_1) d\rho, & x_2-x_1>0,
\end{cases}
\]
From Theorem \ref{thm:CG decomp P} we obtain a product formula for Laguerre polynomials.
\begin{thm} \label{thm:Lag}
The Laguerre polynomials satisfy the following product formula:\\
for $x_1>x_2$
\begin{multline*}
L_{n_1}^{(2k_1-1)}(x_1) L_{n_2}^{(2k_2-1)}(x_2)  = 
\frac{1}{2\pi} \int_0^\infty d_-\, s_{n_2}(\rho^2;k_1-k_2+n_1-n_2+\hf, k_1+k_2-\hf, k_2-k_1+\hf)\\
\times e^{x_2} \ph_{2\rho}^{(2k_2-1,2k_1-1)} \left(\frac{x_2}{x_1-x_2} \right) U(n_2-n_1+k_2-k_1+\hf-i\rho; 1-2i\rho; x_1-x_2) d\rho,
\end{multline*}
for $x_1<x_2$
\begin{multline*}
L_{n_1}^{(2k_1-1)}(x_1) L_{n_2}^{(2k_2-1)}(x_2)  = \frac{1}{2\pi} \int_0^\infty d_+\,s_{n_1}(\rho^2;k_2-k_1+n_2-n_1+\hf, k_1+k_2-\hf, k_1-k_2+\hf)   \\ \times e^{x_1} \ph_{2\rho}^{(2k_1-1,2k_2-1)} \left(\frac{x_1}{x_2-x_1} \right) U(n_1-n_2+k_1-k_2+\hf+i\rho; 1+2i\rho; x_2-x_1)d\rho,
\end{multline*}
for $x_1=x_2=x$
\begin{multline*}
L_{n_1}^{(2k_1-1)}(x) L_{n_2}^{(2k_2-1)}(x)= \frac{1}{\pi} \int_0^\infty s_{n_1}(\rho^2;k_2-k_1+n_2-n_1+\hf, k_1+k_2-\hf, k_1-k_2+\hf)  \\
\times e^{x} \left( d_0\, x^{\hf-k_1-k_2-i\rho} + \overline{d_0} \, x^{\hf-k_1-k_2+i\rho} \right) d\rho,
\end{multline*}
where 
\[
\begin{split}
d_- & =\frac{ (x_1-x_2)^{k_1+k_2-\hf+i\rho} }{n_1!\, n_2!\, \Ga(2k_2)} \left| \frac{ \Ga(k_1+k_2-\hf+i\rho) \Ga(k_2-k_1+n_2-n_1+\hf+i\rho)  } { \Ga(2i\rho) } \right|^2,\\
d_+ & =  \frac{(x_2-x_1)^{k_1+k_2-\hf-i\rho}}{n_1!\, n_2!\, \Ga(2k_1)} \left| \frac{ \Ga(k_1+k_2-\hf+i\rho) \Ga(k_1-k_2+n_1-n_2+\hf+i\rho)  } { \Ga(2i\rho) } \right|^2,\\
d_0 & = \frac{1}{n_1!\, n_2!}\Ga(k_1+k_2-\hf +i\rho) (k_2-k_1+\hf-i\rho)_{n_2-n_1}.
\end{split}
\]
\end{thm}
\begin{proof}
This follows from writing out explicitly $\Ups(e_{n_1} \tensor e_{n_2}) = \Ups_g\big(\int^\oplus S_n(\rho;n_1-n_2) e_{n_1-n_2-L}d\mu^\hf \big)$, $n = \min\{n_1,n_2\}$, where $g$ is given in Theorem \ref{thm:CG decomp P} and using \cite[(3.13)]{GK}
\[
\begin{split}
s_n&(\rho^2;k_1-k_2+\hf, k_1+k_2-\hf, k_2-k_1-p+\hf) = \\
&(-1)^p |(k_1-k_2+\hf+i\rho)_p|^2 s_{n-p} (\rho^2;k_2-k_1+\hf, k_1+k_2-\hf, k_1-k_2+p+\hf).
\end{split}
\]
\end{proof}
\begin{rem}
(i) The confluent hypergeometric $U$-function can be considered as a Whittaker function of the second kind, see \cite[(1.9.6)]{Sl};
\[
W_{k,m}(x) = e^{-\hf x} x^{m+\hf} U(m-k+\hf;2m+1;x).
\]
These Whittaker functions are the kernel in the Whittaker function transform, given by
\[
\begin{cases}
\displaystyle (\mathcal W f )(\la) = \int_0^\infty f(x)  W_{k, i\la}(x) x^{-\frac{3}{2}} dx, \\
\displaystyle f(x) = \frac{1}{2\pi} \int_0^\infty (\mathcal W f )(\la)\, x^{-\hf}  W_{k, i\la}(x) \left| \frac{ \Ga(2i\la) }{ \Ga(\hf-k+i\la) } \right|^2 d\la,
\end{cases}
\]
where $k \leq \hf$. Using the Whittaker function transform we see that the formulas in Theorem \ref{thm:Lag} are a generalization of Koornwinder's formula \cite[(5.14)]{Koo1}, stating that Laguerre polynomials are mapped onto continuous dual Hahn polynomials by the Whittaker function transform. 

(ii) The product formula in Theorem \ref{thm:Lag} has a similar structure as the discontinuous integral for Bessel functions of Weber and Schafheitlin, see \cite[\S13.4]{Wat}.

(iii) Theorem \ref{thm:Lag} can be obtained as a limit case of a bilinear summation formula for Meixner-Pollaczek polynomials, see \cite[Rem.3.2(ii)]{GKR}.
\end{rem}

\end{document}